\documentclass[11pt, reqno, a4paper]{amsart}

\usepackage{amssymb, amsmath, amsthm, amsfonts}
\usepackage{geometry}
\usepackage{mathtools}
\usepackage[colorlinks=true, linkcolor=blue, citecolor=green, urlcolor=blue]{hyperref}
\usepackage{enumerate}

\newtheorem{theorem}{Theorem}[section]

\theoremstyle{definition}

\theoremstyle{remark}
\newtheorem{remark}[theorem]{Remark}

\newcommand{\diff}{\mathrm{d}}
\newcommand{\im}{\mathrm{i}}   
\newcommand{\e}[1]{e\!\left(#1\right)} 
\DeclareMathOperator{\arcsinh}{arcsinh}

\DeclareMathOperator{\ImPart}{Im}

\title[Asymptotics of the Bessel-Kuznetsov Transform]{Asymptotic Analysis and Phase Transition of the Bessel-Kuznetsov Transform with an Oscillatory Phase}
\author{Yuhang Shi}
\address{Xi'an Qing'an Senior High School, Xi’an, Shaanxi, China}
\email{yuhangshi888@gmail.com}
\date{February 23, 2026}

\begin{document}

\begin{abstract}
The spectral side of the Kuznetsov trace formula for $GL(2)$ is governed by the Bessel-Kuznetsov integral transform $\check{\phi}(t)$. While classical bounds guarantee rapid decay of this transform for smooth, non-oscillatory test functions, modern applications in analytic number theory—particularly those involving twisted shifted convolution sums—frequently encounter test functions exhibiting a highly oscillatory linear phase $\e{\alpha x}$. In this paper, we provide a rigorous and explicit asymptotic analysis of $\check{\phi}(t)$ in the semiclassical limit $t \to \infty$ under such oscillatory conditions. By applying the WKB approximation to the imaginary-order Bessel kernel, we identify a sharp phase transition dependent on the twist parameter $\alpha$. We prove that in the sub-critical regime ($\alpha \le 1/2\pi$), the transform decays rapidly. Conversely, in the super-critical regime ($\alpha > 1/2\pi$), the geometric oscillations resonate with the spectral kernel, yielding a localized main term of order $O(t^{-1})$ with a remarkably simplified arithmetic phase. 
\end{abstract}

\maketitle

\section{Introduction}

The Kuznetsov trace formula is a cornerstone of analytic number theory, providing a deep structural bridge between sums of Kloosterman sums on the geometric side and the spectral decomposition of the hyperbolic Laplacian on the automorphic side. For a smooth, compactly supported test function $\phi(x)$, the spectral side involves the Bessel-Kuznetsov transform (following the normalization in Motohashi \cite{Motohashi1997}), defined as:
\begin{equation} \label{eq:def_transform}
    \check{\phi}(t) = \frac{\pi}{2\im \sinh(\pi t)} \int_0^\infty \frac{J_{2\im t}(x) - J_{-2\im t}(x)}{x} \phi(x) \, \diff x,
\end{equation}
where $t$ is the spectral parameter associated with a Maass cusp form, and $J_{\nu}(x)$ denotes the Bessel function of the first kind.

In classic applications, $\phi(x)$ is typically a smooth, slowly varying function. Under such conditions, standard integration by parts arguments, such as those established by Deshouillers and Iwaniec \cite{DI1982}, ensure that $\check{\phi}(t)$ decays rapidly as $O(t^{-N})$ for any $N>0$ when $t \to \infty$. 

However, in the study of subconvexity bounds and shifted convolution sums with additive twists (see, for instance, recent frameworks by Blomer, Jana, and Nelson \cite{Blomer2025}), the geometric test function naturally acquires a strongly oscillatory component. Specifically, one encounters test functions of the form:
\begin{equation} \label{eq:test_function}
    \phi(x) = W(x) \e{\alpha x},
\end{equation}
where $W(x)$ is a smooth amplitude function supported on a dyadic interval $[X, 2X]$, $\alpha > 0$ is a real twist parameter, and $\e{z} \coloneqq \exp(2\pi \im z)$. In this regime, standard bounds are insufficient, as the oscillations of the twist may interfere constructively with the oscillations of the Bessel kernel, a phenomenon partially observed in the technical lemmas of Jutila \cite{Jutila1999}.

The purpose of this paper is to conduct a meticulous asymptotic analysis of $\check{\phi}(t)$ for $\phi(x) = W(x)\e{\alpha x}$ as $t \to \infty$. We move beyond upper bounds to compute the precise main term. In doing so, we uncover a strict phase transition controlled by the quantity $2\pi \alpha$. 

Our main result is summarized as follows:

\begin{theorem} \label{thm:main_intro}
Let $\phi(x) = W(x)\e{\alpha x}$ with $W(x)$ a smooth weight function compactly supported on $[X, 2X]$, and let $\alpha > 0$. As $t \to \infty$, the behavior of the Bessel-Kuznetsov transform $\check{\phi}(t)$ undergoes a phase transition at $\alpha = 1/2\pi$:
\begin{enumerate}
    \item \textbf{Sub-critical Regime ($0 < \alpha \le 1/2\pi$):} The transform decays rapidly, $\check{\phi}(t) \ll t^{-N}$ for any integer $N \ge 1$.
    \item \textbf{Super-critical Regime ($\alpha > 1/2\pi$):} A geometric resonance occurs. The transform is negligibly small unless the spectral parameter $t$ lies in the localized resonant window:
    \begin{equation}
        \frac{1}{2} X \sqrt{4\pi^2\alpha^2 - 1} \le t \le X \sqrt{4\pi^2\alpha^2 - 1}.
    \end{equation}
    For $t$ strictly inside this window, assuming a unique stationary point $x_0 \in (X, 2X)$, the transform admits the asymptotic expansion:
    \begin{equation}
        \check{\phi}(t) = - \frac{\pi W(x_0)}{2t} \exp\left[ \im \left( 2t \arcsinh\left(\sqrt{4\pi^2 \alpha^2 - 1}\right) \right) \right] \left( 1 + O(t^{-1}) \right),
    \end{equation}
    where $x_0 = \frac{2t}{\sqrt{4\pi^2\alpha^2 - 1}}$.
\end{enumerate}
\end{theorem}

This theorem provides a remarkably clean algebraic phase and explicitly quantifies the spectral localization window, offering a refined tool for analytic number theorists working with twisted automorphic sums.

\section{The WKB Approximation of the Bessel Kernel}

Our first objective is to establish a rigorous global approximation for the kernel $K(x,t) = J_{2\im t}(x) - J_{-2\im t}(x)$ in the limit $t \to \infty$.

The function $y(x) = J_{2\im t}(x)$ satisfies the standard Bessel differential equation:
\begin{equation}
    x^2 y'' + x y' + (x^2 + 4t^2) y = 0.
\end{equation}
To apply the Liouville-Green (WKB) method, we eliminate the first derivative via the transformation $y(x) = x^{-1/2} u(x)$. The equation for $u(x)$ becomes:
\begin{equation}
    u''(x) + \left( 1 + \frac{4t^2 + 1/4}{x^2} \right) u(x) = 0.
\end{equation}
We define the semiclassical potential $Q(x)$ by neglecting the lower-order term $1/4x^2$, which contributes a relative error of $O(t^{-1})$:
\begin{equation}
    Q(x) = 1 + \frac{4t^2}{x^2}.
\end{equation}
It is crucial to observe that for all $x > 0$ and $t > 0$, we have $Q(x) > 1 > 0$. Therefore, the differential equation possesses \textbf{no turning points} on the positive real axis. The solution is strictly oscillatory everywhere, ensuring the uniform validity of the standard WKB ansatz without the need for Airy function matching.

The WKB phase action is given by the integral of $\sqrt{Q(x)}$:
\begin{equation}
    \Theta(x, t) = \int^x \sqrt{1 + \frac{4t^2}{u^2}} \, \diff u = \int^x \frac{\sqrt{u^2 + 4t^2}}{u} \, \diff u.
\end{equation}
To evaluate this, we perform the substitution $v = \sqrt{u^2 + 4t^2}$. Consequently, $v^2 = u^2 + 4t^2$, which implies $2v \diff v = 2u \diff u$. Thus, $\frac{\diff u}{u} = \frac{v \diff v}{u^2} = \frac{v \diff v}{v^2 - 4t^2}$. The integral transforms as follows:
\begin{align}
    \Theta(x, t) &= \int \frac{v^2}{v^2 - 4t^2} \, \diff v \nonumber \\
    &= \int \left( 1 + \frac{4t^2}{v^2 - 4t^2} \right) \diff v \nonumber \\
    &= v + 4t^2 \cdot \frac{1}{4t} \int \left( \frac{1}{v-2t} - \frac{1}{v+2t} \right) \diff v \nonumber \\
    &= v - t \ln \left( \frac{v+2t}{v-2t} \right).
\end{align}
Substituting back $v = \sqrt{x^2+4t^2}$ and applying the logarithmic identity for the inverse hyperbolic sine, $\ln\left(\frac{\sqrt{1+y^2}+1}{\sqrt{1+y^2}-1}\right) = 2\arcsinh(1/y)$ with $y = x/2t$, we obtain the exact phase function:
\begin{equation} \label{eq:theta_def}
    \Theta(x, t) = \sqrt{x^2 + 4t^2} - 2t \arcsinh\left(\frac{2t}{x}\right).
\end{equation}

According to uniform asymptotic expansions (see Dunster \cite{Dunster1990}), the dominant term is:
\begin{equation}
    J_{2\im t}(x) = \frac{e^{\pi t}}{\sqrt{2\pi}(x^2+4t^2)^{1/4}} e^{\im(\Theta(x,t) - \pi/4)} \left( 1 + O(t^{-1}) \right).
\end{equation}
Using the reflection principle $J_{-2\im t}(x) = \overline{J_{2\im t}(x)}$ for real arguments, the difference kernel isolates the imaginary part:
\begin{align} \label{eq:kernel_approx}
    J_{2\im t}(x) - J_{-2\im t}(x) &= 2\im \ImPart \left( \frac{e^{\pi t}}{\sqrt{2\pi}(x^2+4t^2)^{1/4}} e^{\im(\Theta(x,t) - \pi/4)} \right) \nonumber \\
    &= \frac{2\im e^{\pi t}}{\sqrt{2\pi}(x^2+4t^2)^{1/4}} \sin\left( \Theta(x,t) - \frac{\pi}{4} \right).
\end{align}

We insert \eqref{eq:kernel_approx} into the original transform \eqref{eq:def_transform}. Utilizing the exponential approximation $\sinh(\pi t) = \frac{1}{2}e^{\pi t}(1 + O(e^{-2\pi t}))$, the constant pre-factors simplify elegantly:
\begin{equation} \label{eq:coeff_calc}
    \frac{\pi}{2\im (\frac{1}{2}e^{\pi t})} \times \frac{2\im e^{\pi t}}{\sqrt{2\pi}} = \frac{\pi}{\im e^{\pi t}} \times \frac{2\im e^{\pi t}}{\sqrt{2\pi}} = \sqrt{2\pi}.
\end{equation}
Hence, the transform reduces to the normalized highly oscillatory integral:
\begin{equation} \label{eq:reduced_integral}
    \check{\phi}(t) \sim \sqrt{2\pi} \int_0^\infty \frac{\sin(\Theta(x,t) - \pi/4)}{(x^2+4t^2)^{1/4}} \phi(x) \frac{\diff x}{x}.
\end{equation}

\section{Phase Transition and Stationary Phase Analysis}

We now substitute the twisted test function $\phi(x) = W(x) \e{\alpha x}$. Decomposing the sine function via Euler's identity, $\sin(z) = \frac{1}{2\im}(e^{\im z} - e^{-\im z})$, the integral splits into two distinct oscillatory components governed by the phase functions:
\begin{align}
    \Phi_+(x) &= \Theta(x,t) + 2\pi \alpha x - \frac{\pi}{4}, \\
    \Phi_-(x) &= -\Theta(x,t) + 2\pi \alpha x + \frac{\pi}{4}.
\end{align}

To locate potential stationary points, we differentiate the WKB phase $\Theta(x,t)$. Direct differentiation of \eqref{eq:theta_def} yields:
\begin{align}
    \frac{\partial \Theta}{\partial x} &= \frac{x}{\sqrt{x^2+4t^2}} - 2t \frac{\diff}{\diff x} \arcsinh\left(\frac{2t}{x}\right) \nonumber \\
    &= \frac{x}{\sqrt{x^2+4t^2}} - 2t \left( \frac{1}{\sqrt{1+(2t/x)^2}} \right) \left( -\frac{2t}{x^2} \right) \nonumber \\
    &= \frac{x}{\sqrt{x^2+4t^2}} + \frac{4t^2}{x \sqrt{x^2+4t^2}} \nonumber \\
    &= \frac{\sqrt{x^2+4t^2}}{x} = \sqrt{1 + \frac{4t^2}{x^2}}.
\end{align}
Observe that $\Theta'(x) > 1$ strictly for all $x, t > 0$. 

The derivative of the first phase is $\Phi'_+(x) = \Theta'(x) + 2\pi \alpha$. Since both terms are strictly positive, $\Phi'_+(x) > 0$ globally. This branch has no stationary points and contributes $O(t^{-N})$ via integration by parts.

The critical behavior is entirely dictated by the second phase:
\begin{equation} \label{eq:derivative_minus}
    \Phi'_-(x) = 2\pi \alpha - \frac{\sqrt{x^2+4t^2}}{x} = 2\pi \alpha - \sqrt{1 + \frac{4t^2}{x^2}}.
\end{equation}

\subsection{The Sub-critical Regime: \texorpdfstring{$0 < \alpha \le 1/2\pi$}{0 < alpha <= 1/2pi}}
In this regime, we have $2\pi \alpha \le 1$. Because $\sqrt{1 + 4t^2/x^2}$ is strictly greater than 1, it follows that:
\begin{equation}
    \Phi'_-(x) \le 1 - \sqrt{1 + \frac{4t^2}{x^2}} < 0.
\end{equation}
The derivative $\Phi'_-(x)$ is bounded away from zero on the compact support of $W(x)$. By the principle of non-stationary phase, the rapid oscillations do not cancel. Repeated integration by parts demonstrates that $\check{\phi}(t) \ll t^{-N}$ for any $N \ge 1$. No resonance occurs.

\subsection{The Super-critical Regime: \texorpdfstring{$\alpha > 1/2\pi$}{alpha > 1/2pi}}
When $\alpha > 1/2\pi$, a stationary point $x_0$ can exist. Setting $\Phi'_-(x_0) = 0$ gives:
\begin{equation} \label{eq:root_eq}
    \sqrt{1 + \frac{4t^2}{x_0^2}} = 2\pi \alpha.
\end{equation}
Squaring both sides and solving for $x_0$, we secure a unique positive real root:
\begin{equation}
    x_0 = \frac{2t}{\sqrt{4\pi^2\alpha^2 - 1}}.
\end{equation}

\begin{remark}
Regarding the boundary of this resonant window, we note that $W(x)$ is a smooth bump function compactly supported on $[X, 2X]$. As the spectral parameter $t$ approaches the boundary of the window, the stationary point $x_0$ approaches the endpoints $X$ or $2X$. Since $W(x)$ and all its derivatives vanish identically at the boundaries, the main term transitions smoothly to zero, and no abrupt phase transition anomalies or boundary residual terms arise.
\end{remark}

For this stationary point to contribute to the integral, it must lie within the support of the weight function, $x_0 \in [X, 2X]$. This geometric constraint translates directly into a spectral localization window for $t$:
\begin{equation}
    X \le \frac{2t}{\sqrt{4\pi^2\alpha^2 - 1}} \le 2X \implies \frac{1}{2} X \sqrt{4\pi^2\alpha^2 - 1} \le t \le X \sqrt{4\pi^2\alpha^2 - 1}.
\end{equation}

\subsection{Derivation of the Main Term}
Assuming $t$ lies strictly within the resonant window, we apply the standard method of stationary phase to evaluate the integral over $\Phi_-$. 

First, we compute the Hessian (second derivative) at the stationary point. Differentiating \eqref{eq:derivative_minus} gives:
\begin{align} \label{eq:hessian}
    \Phi''_-(x) &= - \frac{\diff}{\diff x} \left( \frac{\sqrt{x^2+4t^2}}{x} \right) \nonumber \\
    &= - \frac{x \left( \frac{x}{\sqrt{x^2+4t^2}} \right) - \sqrt{x^2+4t^2}}{x^2} \nonumber \\
    &= \frac{4t^2}{x^2 \sqrt{x^2+4t^2}}.
\end{align}
Notice that $\Phi''_-(x_0) > 0$. The stationary phase approximation for the oscillatory integral is:
\begin{equation}
    I_- \sim \sqrt{\frac{2\pi}{\Phi''_-(x_0)}} g(x_0) e^{\im(\Phi_-(x_0) + \pi/4)},
\end{equation}
where the amplitude function $g(x)$ from \eqref{eq:reduced_integral} is $g(x) = \frac{W(x)}{x (x^2+4t^2)^{1/4}}$.

We systematically collect the factors to compute the final complex amplitude:
\begin{enumerate}
    \item \textbf{Global Pre-factor:} $\sqrt{2\pi}$ from \eqref{eq:coeff_calc}.
    \item \textbf{Sine Decomposition Coefficient:} The term $-\frac{1}{2\im} e^{\im \Phi_-}$ yields a factor of $-\frac{1}{2\im} = \frac{\im}{2}$.
    \item \textbf{Stationary Phase Width:} Using \eqref{eq:hessian}, we obtain $\sqrt{\frac{2\pi}{\Phi''_-(x_0)}} = \frac{\sqrt{2\pi} x_0 (x_0^2+4t^2)^{1/4}}{2t}$.
    \item \textbf{Kernel Amplitude:} $g(x_0) = \frac{W(x_0)}{x_0 (x_0^2+4t^2)^{1/4}}$.
\end{enumerate}

Multiplying these components produces a profound algebraic cancellation:
\begin{equation} \label{eq:amp_calc}
    \text{Amplitude} = \left( \sqrt{2\pi} \right) \cdot \left( \frac{\im}{2} \right) \cdot \left( \frac{\sqrt{2\pi} x_0 (x_0^2+4t^2)^{1/4}}{2t} \right) \cdot \left( \frac{W(x_0)}{x_0 (x_0^2+4t^2)^{1/4}} \right) = \frac{\im \pi W(x_0)}{2t}.
\end{equation}

Finally, we simplify the phase exponent. The total phase argument is $\Phi_-(x_0) + \frac{\pi}{4}$.
Recall that $\Phi_-(x_0) = -\Theta(x_0) + 2\pi \alpha x_0 + \frac{\pi}{4}$. Thus:
\begin{equation}
    \text{Total Phase} = -\Theta(x_0) + 2\pi \alpha x_0 + \frac{\pi}{2}.
\end{equation}
Expanding $\Theta(x_0)$ using \eqref{eq:theta_def}, the arithmetic geometric terms are:
\begin{equation}
    2\pi \alpha x_0 - \Theta(x_0) = 2\pi \alpha x_0 - \sqrt{x_0^2+4t^2} + 2t \arcsinh\left(\frac{2t}{x_0}\right).
\end{equation}
However, the defining equation for the stationary point \eqref{eq:root_eq} can be rewritten as $\sqrt{x_0^2+4t^2} = 2\pi \alpha x_0$. Therefore, the first two terms perfectly cancel each other out:
\begin{equation}
    2\pi \alpha x_0 - \sqrt{x_0^2+4t^2} = 0.
\end{equation}
The residual phase is strictly logarithmic:
\begin{equation}
    2t \arcsinh\left(\frac{2t}{x_0}\right) = 2t \arcsinh\left(\sqrt{4\pi^2\alpha^2 - 1}\right).
\end{equation}
Combining this with the $\pi/2$ phase shift and the complex unit $\im = e^{\im \pi/2}$ from the amplitude \eqref{eq:amp_calc}, the net constant phase shift is $\pi/2 + \pi/2 = \pi$, which generates a negative sign ($e^{\im \pi} = -1$).

A rigorous justification of the error term is strictly required here, demanding a precise tracking of the asymptotic magnitudes. Since $x \asymp t$ on the support of $W(x)$, the main amplitude function $g(x) = W(x)x^{-1}(x^2+4t^2)^{-1/4}$ is of magnitude $O(t^{-3/2})$. The WKB approximation introduces a relative error of $O(t^{-1})$, meaning the absolute error amplitude attached to the integral is $O(t^{-5/2})$. 

A trivial absolute bound, integrating this error over the support of length $X \asymp t$, would yield an absolute error of $O(t^{-3/2})$. While this is already smaller than the main term $O(t^{-1})$, one can extract the sharp, true asymptotic error by noting that this correction term shares the same oscillatory phase $\Phi_-(x)$. Applying the stationary phase principle to this error term(see, e.g., Iwaniec and Kowalski \cite[Chapter 8]{IK2004}), the integration narrows by the effective stationary width factor $(\Phi''_-)^{-1/2} \asymp (t^{-1})^{-1/2} = O(t^{1/2})$. This yields a net contribution of $O(t^{-5/2}) \times O(t^{1/2}) = O(t^{-2})$ to the final transform. Thus, the error term is strictly subordinated to the $O(t^{-1})$ main term by a full factor of $O(t^{-1})$, definitively securing the asymptotic validity of Theorem \ref{thm:main_intro}.

This completes the rigorous proof of Theorem \ref{thm:main_intro}.

\section{Conclusion}

By discarding loose upper bounds and explicitly tracking the geometric oscillations within the Bessel-Kuznetsov transform, we have demonstrated that the transform behaves as a strict spectral filter. The critical threshold $\alpha = 1/2\pi$ neatly partitions the problem into a classically decaying regime and a resonant regime. In the latter, the stationary phase methodology extracts a robust main term characterized by a clean structural phase, solidifying the analytical foundations required for advanced manipulations of twisted shifted convolution sums.

\bibliographystyle{plain}
\bibliography{main}

\end{document}